\documentclass[twoside,11pt]{article}
\usepackage{amssymb}
\usepackage{amsmath}
\usepackage{graphicx}
\usepackage{fancyhdr}

\usepackage{indentfirst}

\newtheorem{theorem}{Theorem}
\newtheorem{corollary}[theorem]{Corollary}
\newtheorem{definition}[theorem]{Definition}
\newtheorem{example}[theorem]{Example}
\newtheorem{proposition}[theorem]{Proposition}
\newtheorem{remark}[theorem]{Remark}
\newenvironment{proof}[1][Proof]{\textbf{#1.} }{\ \rule{0.5em}{0.5em} \bigskip}

\newenvironment{keywords}{\begin{center}
\begin{minipage}[c]{11cm} {\bf Keywords:}} {\end{minipage}
\end{center}}
\newenvironment{msc}{\begin{center}
\begin{minipage}[c]{11cm} {\bf 2000 Mathematics Subject Classification:}} {\end{minipage}
\end{center} \bigskip}
\begin{document}

\title{Nonconservative Noether's Theorem\\ in Optimal Control\footnote{To be
presented at \emph{13th IFAC Workshop on Control Applications of
Optimisation}, 26-28 April 2006, Paris - Cachan, France. Accepted
(19-12-2005) for the Proceedings, IFAC publication, Elsevier Ltd,
Oxford, UK. Research Report CM05/I-54.}}

\author{Gast\~{a}o S. F. Frederico\\
\texttt{gfrederico@mat.ua.pt} \and
Delfim F. M. Torres\\
\texttt{delfim@mat.ua.pt}}
\date{Centre for Research in Optimization and Control\\
Department of Mathematics, University of Aveiro\\
3810-193 Aveiro, Portugal}

\maketitle

\begin{abstract}
We extend Noether's theorem to dynamical optimal control systems
being under the action of nonconservative forces. A systematic way of
calculating conservation laws for nonconservative optimal control
problems is given. As a corollary, the conserved quantities
previously obtained in the literature for nonconservative problems
of mechanics and the calculus of variations are derived.
\end{abstract}

\begin{keywords}
Nonconservative forces, symmetry, conservation laws,
Noether's theorem, optimal control.
\end{keywords}

\begin{msc}
49K05, 49K15, 49S05, 70H33, 37J15, 93C10.
\end{msc}

%%%%%%%%%%%%%%%%%%%%%%%%%%%%%%%%%%%%%%%%%

\section{Introduction}

The concept of symmetry plays an important
role both in Physics and Mathematics.
Symmetries are described by transformations of the system,
which result in the same object after the transformation
is carried out. They are described mathematically
by parameter groups of transformations.
Their importance ranges from fundamental and theoretical aspects
to concrete applications, having profound implications
in the dynamical behavior of the systems, and
in their basic qualitative properties.

Another fundamental notion in Physics and Mathematics
is the one of conservation law.
Typical application of conservation laws
in the calculus of variations and optimal control
is to reduce the number of degrees of freedom,
and thus reducing the problems to a lower dimension,
facilitating the integration of the differential
equations given by the necessary optimality conditions.

Emmy Noether was the first who proved, in 1918, that the
notions of symmetry and conservation law are connected:
when a system exhibits a symmetry, then a conservation law can be obtained.
One of the most important and well known illustrations of this
deep and rich relation, is given by the conservation of energy in Mechanics:
the autonomous Lagrangian $L(q,\dot{q})$,
correspondent to a mechanical system of conservative points,
is invariant under time-translations
(time-homogeneity symmetry), and
\begin{equation}
\label{eq:consEneg}
-L\left(q(t),\dot{q}(t)\right)
+\frac{\partial L}{\partial \dot{q}}\left(q(t),\dot{q}(t)\right) \cdot \dot{q}(t)
\equiv \text{constant}
\end{equation}
follows from Noether's theorem,
\textrm{i.e.}, the total energy of a conservative
closed system always remain constant in time, ``it cannot be created
or destroyed, but only transferred from one form into another''.
Expression \eqref{eq:consEneg}
is valid along all the Euler-Lagrange extremals $q(\cdot)$
of an autonomous problem of the calculus of variations.
The conservation law \eqref{eq:consEneg} is known
in the calculus of variations as the 2nd Erdmann necessary condition;
in concrete applications, it gains different interpretations:
conservation of energy in Mechanics;
income-wealth law in Economics;
first law of Thermodynamics; etc.
The literature on Noether's theorem is vast,
and many extensions of the classical results of Emmy Noether
are now available for the more general
setting of optimal control (see \textrm{e.g.}
\cite{CD:Djukic:1972,Gogodze88,CD:JMS:Torres:2002,CD:PortMath:Torres:2004},
and references therein). Here we remark that in all those results
conservation laws always refer to closed systems.

It turns out that in practical terms closed systems do not exist:
forces that do not store energy, so-called nonconservative
or dissipative forces, are always present in real systems.
Friction is a nonconservative force, but others do exist.
Any friction-type force, like air resistance,
is a nonconservative force. Nonconservative forces remove energy from the systems
and, as a consequence, the conservation law \eqref{eq:consEneg} is broken.
This explains, for instance, why the innumerable
``perpetual motion machines'' that have been proposed fail.
In presence of external nonconservative forces,
Noether's theorem and respective conservation laws cease to be valid.
However, it is still possible to obtain a Noether-type theorem
which covers both conservative (closed system)
and nonconservative cases \cite{CD:Djukic:1980,CD:LiQun:2003}.
Roughly speaking, one can prove
that Noether's conservation laws are still valid if a new term,
involving the nonconservative forces, is added to the standard
conservation laws.

Here we extend previous nonconservative results
\cite{CD:Djukic:1980,CD:LiQun:2003}
to the wider context of optimal control.
For that, and differently from \cite{CD:Djukic:1980,CD:LiQun:2003},
where the Lagrangian formalism is considered,
we adopt an Hamiltonian point of view.

%%%%%%%%%%%%%%%%%%%%%%%%%%%%%%%%%%%%%%%%%%%%%%%%%%

\section{Preliminaries}

Let us consider the optimal control problem in Lagrange form:
\begin{gather}
\label{eq:JO}
I[q(\cdot),u(\cdot)] =\int_a^b L\left(t,q(t),u(t)\right) dt \longrightarrow \min \, , \tag{P} \\
\dot{q}(t)=\varphi\left(t,q(t),u(t)\right) \, , \notag
\end{gather}
together with some boundary conditions on $q(\cdot)$.
In problem \eqref{eq:JO} $\dot{q}=\frac{dq}{dt}$, and the Lagrangian
$L : [a,b]\times \mathbb{R}^{n} \times \mathbb{R}^{m} \rightarrow \mathbb{R}$
and the velocity vector
$\varphi : [a,b] \times \mathbb{R}^{n} \times \mathbb{R}^{m} \rightarrow\mathbb{R}^n$
are assumed to be $C^{1}$ functions with respect to all the arguments.
In agreement with the calculus of variations,
we assume the admissible state trajectories
to be piecewise smooth, and the admissible control functions
to be piecewise constant with no restrictions on their values:
$q(\cdot) \in PC^1([a,b];\mathbb{R}^n)$, $u(\cdot) \in PC([a,b];\mathbb{R}^m)$.

\begin{remark}
There is no classical Hamiltonian theory for constrained
variational problems. For an ongoing attempt to develop it,
we refer the reader to \cite[Ch.~6]{livroAMS2005Clarke}.
\end{remark}

\begin{remark}
The fundamental problem of the calculus of variations,
\begin{equation}
\label{eq:pfcv}
I[q(\cdot)] = \int_a^b L\left(t,q(t),\dot{q}(t)\right) \longrightarrow \min \, ,
\end{equation}
is a particular case of problem \eqref{eq:JO}: in that case $\varphi(t,q,u)=u$.
The problems of the calculus of variations with higher-order
derivatives are also easily written in the optimal control form \eqref{eq:JO}.
For example, the problem of the calculus of variations with derivatives
of second order,
\begin{equation}
\label{eq:pcvo2}
I[q(\cdot)] = \int_a^b L\left(t,q(t),\dot{q}(t),\ddot{q}(t)\right) \longrightarrow \min \, ,
\end{equation}
is equivalent to problem
\begin{gather*}
I[q^0(\cdot),q^1(\cdot),u(\cdot)]
= \int_a^b L\left(t,q^0(t),q^1(t),u(t)\right) \longrightarrow \min \, , \\
\begin{cases}
\dot{q}^0(t) = q^1(t) \, , \\
\dot{q}^1(t) = u(t) \, .
\end{cases}
\end{gather*}
\end{remark}

In the fifties of the twentieth century,
L.S.~Pontryagin and his collaborators proved the main
necessary optimality condition for optimal control problems:
the famous Pontryagin Maximum Principle \cite{CD:MR29:3316b}.

\begin{definition}[Process]
An admissible pair $(q(\cdot),u(\cdot))$ satisfying the control system
$\dot{q}(t)=\varphi\left(t,q(t),u(t)\right)$ of \eqref{eq:JO},
$t \in [a,b]$, is called a \emph{process}.
\end{definition}

\begin{theorem}[Pontryagin Maximum Principle-- \cite{CD:MR29:3316b}]
\label{th:P}
If $(q(\cdot),u(\cdot))$ is an optimal process for problem
\eqref{eq:JO}, then there exists a vectorial function
$p(\cdot) \in PC^1([a,b];\mathbb{R}^n)$ such that the following conditions hold:
\begin{itemize}
\item the Hamiltonian system
\begin{equation}
\label{eq:Ham}
\begin{cases}
\dot{q}(t)&=\frac{ \partial{{\cal H}} }{\partial {p}}(t, q(t), u(t),p(t)) \, , \\
\dot{p}(t) &= -\frac{ \partial{{\cal H}} }{\partial {q}}(t, q(t), u(t), p(t)) \, ;
\end{cases}
\end{equation}
\item the stationary condition
\begin{equation}
\label{eq:CE}
 \frac{ \partial{{\cal H}} }{\partial {u}}(t, q(t), u(t), p(t))=0 \, ;
\end{equation}
\end{itemize}
with the Hamiltonian ${\cal H}$ defined by
\begin{equation}
\label{eq:H}
{\cal H}\left(t,q,u,p\right)
=-L\left(t,q,u\right)+p \cdot \varphi\left(t,q,u\right) \, .
\end{equation}
\end{theorem}

\begin{remark}
In mechanics, $p$ corresponds to the \emph{generalized momentum}.
In the language of optimal control $p$
is called the \emph{adjoint variable}.
\end{remark}

\begin{definition}
\label{def:extPont}
Any triplet $(q(\cdot),u(\cdot),p(\cdot))$ satisfying
the conditions of Theorem~\ref{th:P} is called
a \emph{Pontryagin extremal}.
\end{definition}

\begin{remark}
The Pontryagin Maximum Principle is more general than
we state it here. Written as in Theorem~\ref{th:P},
the Pontryagin Maximum Principle is also known
as \emph{Hestenes Theorem}. In particular,
we are only considering normal Pontryagin extremals.
\end{remark}
From the Hamiltonian system \eqref{eq:Ham}
and the stationary condition \eqref{eq:CE},
it follows that
\begin{equation}
\label{eq:H2} \frac{d{\cal H}}{dt}(t, q(t), u(t), p(t))=\frac{
\partial{{\cal H}} }{\partial {t}}(t, q(t), u(t), p(t)) \, .
\end{equation}
When the optimal control problem \eqref{eq:JO}
is autonomous (when the Hamiltonian
\eqref{eq:H} does not depend explicitly on time $t$) one
obtains from \eqref{eq:H2} the conservation law
\begin{equation}
\label{eq:Hc}
{\cal H}(t, q(t), u(t), p(t)) = const \, .
\end{equation}

For the fundamental problem of the calculus of variations \eqref{eq:pfcv}
$(\varphi=u \Rightarrow {\cal H} = -L + p \cdot u)$
one obtains from the Pontryagin Maximum Principle:
\begin{gather*}
 \dot{q}=\frac{ \partial{{\cal H}} }{\partial {p}}=u \, ,\\
 \dot{p} = -\frac{ \partial{{\cal H}} }{\partial {q}}=\frac{\partial L}{\partial q} \, ,\\
 \frac{ \partial{{\cal H}}}{\partial {u}}=0 \Leftrightarrow
 p=\frac{\partial L}{\partial u}\Rightarrow\dot{p}=\frac{d}{dt}\frac{\partial L}{\partial u} \, .
\end{gather*}
Comparing the two expressions for $\dot{p}$,
one arrives to the Euler-Lagrange differential equations:
\begin{equation}
\label{eq:EL}
 \frac{d}{{dt}}\frac{{\partial L}} {{\partial u}} =
\frac{{\partial L}} {{\partial q}}.
\end{equation}
When the fundamental problem of the calculus of variations
is autonomous, equality \eqref{eq:Hc} reduces to \eqref{eq:consEneg}.

In the presence of nonconservative forces $Q(t,q(t),\dot{q}(t))$,
\textrm{i.e.} forces which are not equivalent to the gradient of a potential,
like friction and drag, the Euler-Lagrange equations \eqref{eq:EL}
are no longer valid, and it is well-known that they must be substituted by
\begin{equation}
\label{eq:ELNC}
 \frac{d}{{dt}}\frac{{\partial L}} {{\partial u}} =
\frac{{\partial L}} {{\partial q}}+Q
\end{equation}
(see \textrm{e.g.} \cite{CD:LiQun:2003}). In Physics the Hamiltonian formalism
is not common for the nonconservative case. However, it is clear that
Theorem~\ref{th:P} and property \eqref{eq:H2} must be also changed
when considering the influence of external dissipative forces $Q$.
This will be addressed in \S\ref{sec:mr}.

Following \cite{CD:Djukic:1972},
the notion of invariance for problem \eqref{eq:JO}
is defined in terms of the Hamiltonian, by introducing
the augmented functional
 \begin{equation*}
\label{eq:J} J[q(\cdot),u(\cdot),p(\cdot)]
= \int_a^b \left[{\cal H}\left(t,q(t),u(t),p(t)\right)-p(t)\cdot\dot{q}(t)\right]dt \, ,
\end{equation*}
where ${\cal H}$ is given by \eqref{eq:H}.

\begin{definition}[Invariance up to a gauge term -- \textrm{cf.} \cite{CD:Djukic:1972}]
\label{def:inv:gt}
An optimal control problem \eqref{eq:JO} is said to be invariant under
the $\varepsilon$-parameter local group of transformations
\begin{equation}
\label{eq:trf:inf}
\begin{cases}
\bar{t}(t) = t+\varepsilon\tau(t, q(t), u(t), p(t)) + o(\varepsilon) \, , \\
\bar{q}(t) = q(t)+\varepsilon\xi(t, q(t), u(t), p(t)) + o(\varepsilon) \, , \\
\bar{u}(t) = u(t)+\varepsilon\sigma(t, q(t), u(t), p(t)) + o(\varepsilon) \, , \\
\bar{p}(t) = p(t)+\varepsilon\alpha(t, q(t), u(t), p(t))+ o(\varepsilon) \, , \\
\end{cases}
\end{equation}
if, and only if, there exists a function $\Lambda$ such that
\begin{equation}
\label{eq:H3}
\left[{\cal H}(\bar{t},\bar{q},\bar{u},\bar{p})-\bar{p} \cdot
\dot{\bar{q}}\right]d\bar{t}
=\left[{\cal H}(t,q,u,p)-p \cdot \dot{q}\right]dt
+ \varepsilon d\Lambda(t,q,u,p) \, .
\end{equation}
\end{definition}

\begin{remark}
Function $\Lambda$ of Definition~\ref{def:inv:gt}
is called a \emph{gauge term} in the Physics literature.
In the particular case $\Lambda=0$, one obtains
the concept of \emph{absolute invariance}.
\end{remark}

\begin{remark}
We can write equation \eqref{eq:H3} in the following way:
\begin{equation}
\label{eq:H5}
[{\cal H}(\bar{t},\bar{q},\bar{u},\bar{p})-\bar{p}\cdot \dot{\bar{q}}]
\frac{d\bar{t}}{dt}
=\left[{\cal H}(t,q,u,p)-p\cdot \dot{q}\right]+\varepsilon \frac{d \Lambda}{dt}.
\end{equation}
\end{remark}

Functions $\tau$, $\xi$, $\sigma$, and $\alpha$
are known as the \emph{infinitesimal generators}
of the invariance-transformations \eqref{eq:trf:inf}.
Next theorem asserts that generators
are sufficient to define invariance.

\begin{theorem}[Necessary and sufficient condition of invariance]
\label{th:cnsi}
Problem \eqref{eq:JO} is said to be invariant
up to the gauge term $\Lambda$ if, and only if,
\begin{equation}
\label{eq:H4}
\tau\frac{\partial{\cal H}}{\partial
t}+\xi\cdot\frac{\partial{\cal H}}{\partial
q}+\sigma\cdot\frac{\partial{\cal H}}{\partial
u}+\alpha\cdot(\frac{\partial{\cal H}}{\partial
p}-\dot{{q}})-\dot{\xi}\cdot p+\dot{\tau}{\cal H}
= \frac{d \Lambda}{dt}.
\end{equation}
\end{theorem}
\begin{proof}
Having in mind that for $\varepsilon=0$ one has $\bar{t} = t$,
$\bar{q} = q$, $\bar{u} = u$, $\bar{p} = p$,
we differentiate \eqref{eq:H5} with respect to
$\varepsilon$ and then put $\varepsilon=0$.
\end{proof}

\begin{definition}
A quadruple $\left(\tau,\xi,\sigma,\alpha\right)$
is said to be a \emph{symmetry} of the optimal control problem
\eqref{eq:JO} if it satisfies condition \eqref{eq:H4}
for a certain $\Lambda$. We talk about \emph{exact symmetries}
if $\Lambda = 0$.
\end{definition}

\begin{remark}
It is possible to use a modern computer algebra system
to compute the symmetries of an optimal control problem
\eqref{eq:JO} in an automatic way -- see \cite{GouveiaTorresCMAM}.
\end{remark}

\begin{remark}
For the fundamental problem of the calculus of variations \eqref{eq:pfcv},
the necessary and sufficient condition of invariance \eqref{eq:H4} takes
the well-known form (\textrm{cf.} \textrm{e.g.} \cite[pp.~429]{CD:Logan:1987})
\begin{equation*}
\tau \frac{\partial L}{\partial t} + \xi \cdot \frac{\partial L}{\partial q}
+ \frac{\partial L}{\partial \dot{q}} \cdot \left(\dot{\xi} - \dot{q} \dot{\tau}\right)
+ \dot{\tau} L = \frac{d \Lambda}{dt} \, .
\end{equation*}
\end{remark}

\begin{remark}
For the problem of the calculus of variations with derivatives
of second order \eqref{eq:pcvo2}, the necessary and sufficient
condition of invariance \eqref{eq:H4} takes the form
(\textrm{cf.} \cite[Lemma~5.5]{CD:JMS:Torres:2004})
\begin{equation*}
\tau \frac{\partial L}{\partial t} + \xi_0 \cdot \frac{\partial L}{\partial q}
+ \xi_1 \cdot \frac{\partial L}{\partial \dot{q}}
+ \left(\frac{\partial L}{\partial \dot{q}}
-\frac{d}{dt}\frac{\partial L}{\partial \ddot{q}}\right) \cdot \left(\dot{\xi}_0 - \dot{\tau}\dot{q}\right)
+ \frac{\partial L}{\partial \ddot{q}} \cdot \left(\dot{\xi}_1 - \dot{\tau}\ddot{q}\right)
+ \dot{\tau} L = \frac{d \Lambda}{dt} \, .
\end{equation*}
\end{remark}

%%%%%%%%%%%%%%%%%%%%%%%%%%%%%%%%%%%%%%%%%

\section{Main Result}
\label{sec:mr}

We begin by introducing the notion
of nonconservative Hamiltonian system.
Such concept must lead us to equations \eqref{eq:ELNC}
in the particular case of the fundamental problem
of the calculus of variations under presence
of nonconservative forces $Q$.

\begin{definition}[Nonconservative Hamiltonian System]
\label{def:SHNC}
We define the \emph{nonconservative Hamiltonian system} by
\begin{equation}
\label{eq:HNC}
\begin{cases}
\dot{q}(t)&=\frac{ \partial{{\cal H}} }{\partial {p}}(t, q(t), u(t),p(t)) \, ,\\
\dot{p}(t) &= -\frac{ \partial{{\cal H}} }{\partial {q}}(t,q(t),u(t), p(t))+Q(t,q(t),u(t)) \, ,
\end{cases}
\end{equation}
where ${\cal H}$ is given as in \eqref{eq:H}.
\end{definition}

\begin{remark}
For the particular case of the fundamental
problem of the calculus of variations ($\varphi=u$),
the nonconservative Hamiltonian system \eqref{eq:HNC},
together with the stationary condition \eqref{eq:CE},
lead us to the nonconservative Euler-Lagrange equations \eqref{eq:ELNC}:
\begin{gather*}
 \dot{q}=\frac{ \partial{{\cal H}} }{\partial {p}}=u \, ,\\
 \dot{p} = -\frac{ \partial{{\cal H}} }{\partial
{q}}+Q=\frac{\partial L}{\partial q}+Q \, ,\\
 \frac{ \partial{{\cal H}}
}{\partial {u}}=0 \Leftrightarrow p=\frac{\partial L}{\partial
u}\Rightarrow\dot{p}=\frac{d}{dt}\frac{\partial L}{\partial u} \, ,
\end{gather*}
and comparing the two expressions for $\dot{p}$
we obtain equation \eqref{eq:ELNC}.
\end{remark}

Similarly to Definition~\ref{def:extPont},
we introduce now the notion of nonconservative extremal.

\begin{definition}[Nonconservative Extremal]
Any triplet $(q(\cdot),u(\cdot),p(\cdot))$, satisfying
the stationary condition \eqref{eq:CE} and the
nonconservative Hamiltonian system \eqref{eq:HNC},
will be called a \emph{nonconservative extremal}.
\end{definition}

\begin{proposition}
\label{prp:tdH}
The following property holds along
the nonconservative extremals:
\begin{multline}
\label{eq:H8}
\frac{d{\cal H}}{dt}(t, q(t), u(t), p(t)) \\
= \frac{\partial{{\cal H}} }{\partial {t}}(t, q(t), u(t),
p(t))+Q(t,q(t),u(t))\cdot\frac{ \partial{{\cal H}} }{\partial
{p}}(t, q(t), u(t),p(t)) \, .
\end{multline}
\end{proposition}
\begin{proof}
Equality \eqref{eq:H8} is a simple consequence of the
the stationary condition \eqref{eq:CE} and the
nonconservative Hamiltonian system \eqref{eq:HNC}.
\end{proof}

\begin{remark}
In the particular case $Q = 0$ (conservative case),
the nonconservative Hamiltonian system \eqref{eq:HNC}
takes the form \eqref{eq:Ham}, and the
set of nonconservative extremals coincide
with the set of Pontryagin extremals. In that situation,
property \eqref{eq:H8} equals \eqref{eq:H2}.
\end{remark}

\begin{definition}[Nonconservative Constants of Motion]
\label{def:NCLs}
We say that a function $C(t,q,u,p)$ is
a \emph{nonconservative constant of motion}
if it is preserved along any nonconservative extremal
$(q(\cdot),u(\cdot),p(\cdot))$:
\begin{equation}
\label{eq:CL:NC}
C(t,q(t),u(t),p(t)) = c \, , \quad c \text{ constant } \, ,
\quad \forall \; t \in [a,b] \, .
\end{equation}
Equation \eqref{eq:CL:NC} is then said to be
a \emph{nonconservative conservation law}.
 \end{definition}

Generalizations of the classical Noether's theorem
of the calculus of variations include: (i) (conservative)
Noether-type theorems for the more general framework of optimal control,
asserting the existence of a preserved quantity
along the Pontryagin extremals, whenever a symmetry occurs -- see
\cite{CD:Djukic:1972,Gogodze88,CD:JMS:Torres:2002,CD:PortMath:Torres:2004};
(ii) Noether-type theorems for the nonconservative
calculus of variations, asserting that existence
of a symmetry implies the existence of a preserved
quantity along the solutions of the nonconservative Euler-Lagrange equations
\eqref{eq:ELNC} -- see \cite{CD:Djukic:1980,CD:LiQun:2003}.
Next theorem extends the mentioned Noether-type results:
to each symmetry of the optimal control problem \eqref{eq:JO}
there exists a nonconservative conservation law in the
sense of Definition~\ref{def:NCLs}. Essentially, for $Q = 0$ one gets
from Theorem~\ref{th:TNNC} the results obtained in
\cite{CD:Djukic:1972,Gogodze88,CD:JMS:Torres:2002,CD:PortMath:Torres:2004};
restricting ourselves to the problems of the calculus of variations,
one gets the results found in \cite{CD:Djukic:1980,CD:LiQun:2003}.

\begin{theorem}[Nonconservative Noether's Theorem]
\label{th:TNNC}
If $\left(\tau,\xi,\sigma,\alpha\right)$
is a symmetry of the optimal control problem
\eqref{eq:JO}, and there exists a function $f=f(t,q,u)$ such that
\begin{equation}
\label{eq:H6}
\frac{df}{dt}=Q\cdot(\xi-\tau \dot{q}) \, ,
\end{equation}
where $Q$ denotes the external nonconservative forces
acting on the system, then
\begin{equation}
\label{eq:H7}
C(t,q,u,p)={\cal H}(t,q,u,p)\tau-p \cdot \xi+
f(t,q,u)-\Lambda(t,q,u,p)
\end{equation}
is a nonconservative constant of motion.
\end{theorem}
\begin{remark}
Similarly to \cite{CD:Djukic:1972,Gogodze88,CD:JMS:Torres:2002,CD:PortMath:Torres:2004},
only the generators $\tau$ and $\xi$, corresponding to the transformations
of the time and state variables, appear in Noether's conservation laws
(\textrm{cf.} expression \eqref{eq:H7}).
\end{remark}
\begin{proof} Substituting (\textrm{cf.} Proposition~\ref{prp:tdH})
\begin{equation*}
\tau\frac{\partial{\cal H}}{\partial t}+\dot{\tau}{\cal
H}-p\cdot\dot{\xi}=\frac{d}{dt}({\cal H}\tau-p\cdot\xi)-\tau
Q\cdot\frac{\partial{\cal H}}{\partial p}+\dot{p}\cdot\xi
\end{equation*}
into \eqref{eq:H4}, and using conditions \eqref{eq:CE},
\eqref{eq:HNC}, and \eqref{eq:H6}, we obtain successively:
\begin{equation*}
\label{eq:H10}
 0=\frac{d}{dt}\left({\cal
H}\tau-p\cdot\xi-\Lambda\right)-\tau Q\cdot\frac{\partial{\cal H}}{\partial
p}+\left(\frac{\partial{\cal H}}{\partial
q}+\dot{p}\right)\cdot\xi \, ,
\end{equation*}
\begin{equation*}
0=\frac{d}{dt}({\cal
H}\tau-p\cdot\xi-\Lambda)+\xi\cdot(-\dot{p}+Q+\dot{p})-\tau\dot{q}\cdot
Q \, ,
\end{equation*}
\begin{equation*}
0=\frac{d}{dt}({\cal H}\tau-p\cdot\xi-\Lambda)+
Q\cdot(\xi-\tau\dot{q})=\frac{d}{dt}({\cal
H}\tau-p\cdot\xi+f-\Lambda) \, .
\end{equation*}
\end{proof}

\begin{remark}
Proof of Theorem~\ref{th:TNNC} is very simple.
This shows, in the opinion of the authors,
that the Hamiltonian formalism provides the
natural language to Noether's theory. It is strange,
from the mathematical point of view, why the Hamiltonian
approach is not used in Physics with respect
to nonconservative systems.
\end{remark}
As corollaries, we obtain the previous results known in the literature.
\begin{corollary}\textbf{(Optimal Control version of Noether's Theorem
-- \textrm{cf. e.g.} \cite{CD:Djukic:1972,CD:JMS:Torres:2002})}
In the absence of nonconservative forces (\textrm{i.e.} in the
conservative case $Q = 0$), if the optimal control problem
\eqref{eq:JO} is invariant under the one-parameter transformations
\eqref{eq:trf:inf}, then function
\begin{equation}
\label{eq:H9}
C(t,q,u,p)={\cal H}\tau-p\cdot \xi-\Lambda
\end{equation}
is preserved along any Pontryagin extremal.
\end{corollary}
\begin{proof}
When $Q = 0$, \eqref{eq:H6} implies that $f$ is a constant,
and the conservation law associated with the constant of motion
\eqref{eq:H7} is equivalent to the one associated with \eqref{eq:H9}.
\end{proof}

\begin{corollary}\textbf{(Nonconservative Noether's theorem of the calculus
of variations -- \textrm{cf.} \cite{CD:LiQun:2003})}
\label{cor:H11}
For the fundamental problem of the calculus of variations
\eqref{eq:pfcv} the constant of motion \eqref{eq:H7} is equivalent to
\begin{equation}
\label{eq:H11}
C(t,q,\dot{q})=\tau \left( L-\dot{q}\cdot\frac{\partial
L}{\partial \dot{q}}\right)+\frac{\partial L}{\partial
\dot{q}}\cdot\xi-f+\Lambda \, ,
\end{equation}
that is, under the invariance hypotheses of Theorem~\ref{th:TNNC},
expression \eqref{eq:H11}
is preserved along all the solutions $q(\cdot)$
of the nonconservative Euler-Lagrange equations \eqref{eq:ELNC}.
\end{corollary}
\begin{proof} For the fundamental problem of the calculus
of variations, $\varphi = u$ and
$p = \frac{\partial L}{\partial \dot{q}}$, so that the Hamiltonian
takes the form ${\cal H}=-L+\frac{\partial
L}{\partial \dot{q}}\cdot\dot{q}$. Substituting this expression in
\eqref{eq:H7}, we obtain the desired conclusion.
\end{proof}

 \begin{corollary}\textbf{(Nonconservative Noether's theorem
 for problems of the calculus of variations with second-order derivatives
-- \textrm{cf.} \cite{CD:Djukic:1980})}
For the higher-order problem of the calculus of variations
\eqref{eq:pcvo2}, the constant of motion \eqref{eq:H7} is equivalent to
\begin{equation}
\label{eq:H14}
L\tau+\left(
\frac{\partial L}{\partial \dot{q}}-\frac{d}{dt}\frac{\partial
L}{\partial \ddot{q}}\right) \cdot (\xi_{0}-\dot{q}\tau)+\frac{\partial
L}{\partial \ddot{q}} \cdot (\xi_{1}-\ddot{q}\tau)-f + \Lambda \, .
\end{equation}
 \end{corollary}
\begin{proof}
For the problem of the calculus of variations with second-order derivatives,
one has ${\cal H}\left(t,q^{0},q^{1},u,p^{0},p^{1}\right)=
-L(t,q^{0},q^{1},u)+p^{0} q^{1}+p^{1}u$,
$q^{0}(t)=q(t)$, $q^{1}(t)=\dot{q}(t)$, $u(t)=\ddot{q}(t)$.
Using these equalities, it follows from the Pontryagin Maximum Principle that
\begin{gather*}
\frac{\partial {\cal H}}{\partial u}=0\Leftrightarrow
{p}^{1}=\frac{\partial L}{\partial \ddot{q}} \, ,\\
\dot{p}^{0}=-\frac{\partial {\cal H}}{\partial
q^{0}}=\frac{\partial
L}{\partial q} \, ,\\
\dot{p}^{1}=-\frac{\partial {\cal H}}{\partial
q^{1}} \Leftrightarrow
p^{0}=\frac{\partial L}{\partial
\dot{q}}-\frac{d}{dt}\frac{\partial L}{\partial
\ddot{q}} \, .
\end{gather*}
In this case, the constant of motion \eqref{eq:H7} takes the form
\begin{equation*}
C={\cal H}\tau-p^{0} \cdot \xi_{0}-p^{1} \cdot \xi_{1}+f - \Lambda \, ,
\end{equation*}
and substituting ${\cal H}$, $p^{0}$ and $p^{1}$
by its expressions, the intended result is obtained.
\end{proof}

%%%%%%%%%%%%%%%%%%%%%%%%%%%%%%%%%%%%%%%%%

\section{Examples}

We now illustrate the application of our result
to some concrete problems of nonconservative classical mechanics.
In all the examples, one can use the computational tools
\cite{GouveiaTorresCMAM} to determine the symmetries.

\begin{example}(Forced Oscillations -- \textrm{cf.} \cite[pp.~114--115]{CD:Maia:2000}).
\label{ex:Maia}
Let us consider the problem of vertical oscillations
of a body of mass $m$, connected to a spring of ignorable mass
with constant of elasticity $k$, under the action
of a sinusoidal nonconservative force
$Q(t)=Fe^{i\omega t}$, where $F$ and $\omega$ are positive constants.
In this situation, the associated variational problem
is given by (\textrm{cf.} \cite[pp.~114--115]{CD:Maia:2000})
 \begin{gather*}
I[q(\cdot),u(\cdot)] =\frac{1}{2}\int_0^T (mu^2-kq^2) dt \longrightarrow \min \, , \\
\dot{q}(t)=u(t) \, .
\end{gather*}
The problem exhibits the following exact symmetries:
$\left(\tau,\xi,\sigma,\alpha\right) = \left(c,0,0,0\right)$,
where $c$ is an arbitrary constant. From Theorem~\ref{th:TNNC}
we conclude that
\begin{equation}
\label{eq:OC:H12}
-\frac{1}{2} \left( m u(t)^2 -  k q(t)^2 \right) + p(t) u(t)
- \int \dot{q}(t) F \mathrm{e}^{i w t} dt = const \, , \quad
t \in [0,T] \, ,
\end{equation}
is a nonconservative conservation law. Corollary~\ref{cor:H11}
permits to write the conclusion \eqref{eq:OC:H12} in the notation
of the calculus of variations:
\begin{equation*}
\label{eq:H12}
C(t,q,\dot{q})=\frac{1}{2}(m\dot{q}^2+kq^2)-\int
F\dot{q}e^{i\omega t}dt
\end{equation*}
is a nonconservative constant of motion
(nonconservative extremals are defined by
$u(t) = \dot{q}(t)$, $\dot{p}(t) = F \mathrm{e}^{i w t} - k q(t)$,
$p(t) = m u(t)$). In this simple example, such fact
is easily verified by direct application of Definition~\ref{def:NCLs}:
$\frac{dC}{dt}=0 \Leftrightarrow
\dot{q} \left( m\ddot{q}+ kq - Fe^{i\omega t} \right) = 0$,
which is a truism along the nonconservative Euler-Lagrange
extremals -- the solutions of the equation \eqref{eq:ELNC}
$m\ddot{q}+ kq - Fe^{i\omega t} = 0$.
\end{example}

\begin{example}[\cite{CD:LiQun:2003}]
We apply our results to the example studied in \cite{CD:LiQun:2003}:
a dynamical system with Lagrangian
$L(t,q,\dot{q}) = \dot{q}^2/2$,
subjected to the nonconservative force $Q = \dot{q}^2$.
The associated variational problem is given by
 \begin{gather*}
I[q(\cdot),u(\cdot)] =\frac{1}{2}\int_0^T u^2 dt \longrightarrow \min \, , \\
\dot{q}(t)=u \, ,
\end{gather*}
and from Theorem~\ref{th:cnsi} we get the exact symmetries
\begin{equation*}
\left(\tau,\xi,\sigma,\alpha\right)
= \left(2 c_1 t + c_2,c_1 q + c_3,- c_1 u,- c_1 p\right) \, ,
\end{equation*}
where $c_i$, $i = 1,2,3$, are constants. From Theorem~\ref{th:TNNC},
we obtain the nonconservative conservation law
\begin{multline}
\label{eq:lcgex2}
\left(c_1 q(t)+c_3\right) p(t)
+\left(\frac{1}{2} u(t)^2-p(t) u(t)\right) \left(2 c_1 t+c_2\right)\\
+ \int \left((2 c_1 t+c_2) \dot{q}(t)-c_1 q(t)-c_3\right) u(t)^2 dt = const \, .
\end{multline}
The equation of motion is defined by
$u(t) = \dot{q}(t)$, $\dot{p}(t) = u(t)^2$, $p(t) = u(t)$
(which is equivalent to the nonconservative Euler-Lagrange
equation \eqref{eq:ELNC}: $\ddot{q}(t) = \dot{q}^2$), and
we obtain the nonconservative conservation law \eqref{eq:lcgex2}
in terms of the calculus of variations \eqref{eq:H11}:
\begin{multline*}
c_1 \left[ q(t) \dot{q}(t) - t \dot{q}(t)^{2} \right]
-\frac{1}{2} \, c_2 \dot{q}(t)^{2} + c_3 \dot{q}(t) \\
 +\int \! \left[ 2 c_1 t \dot{q}(t) - c_1 q(t)
+ c_2 \dot{q}(t) - c_3 \right] \dot{q}(t)^{2} dt = const \, .
\end{multline*}
Similarly to Example~\ref{ex:Maia}, also here it is possible
to verify the validity of \eqref{eq:lcgex2} by definition:
the nonconservative extremals $(q(\cdot),u(\cdot),p(\cdot))$ are
given by $q(t) = k_1 - \ln(t-k_2)$, $u(t) = \frac{1}{k_2-t}$,
$p(t) = \frac{1}{k_2-t}$, where $k_1$ and $k_2$ are constants to be determined
from the boundary conditions, and substituting the extremals
in \eqref{eq:lcgex2}, we obtain the tautology $0 = const$.
\end{example}

To finish the illustration of our methods,
we consider a generalized mechanical system with
one degree of freedom, whose Lagrangian and
nonconservative force depends on higher-order derivatives.
\begin{example}[\cite{CD:Djukic:1980}]
The following problem is borrowed from \cite[\S 4]{CD:Djukic:1980}:
\begin{gather*}
L=\frac{1}{2} \left( \ddot{q}^2+ a\dot{q}^2 + bq^2 \right) \, , \\
Q=\mu\dot{q}+\left(\frac{\mu}{a}\right)^2\ddot{q}-2 \left(\frac{\mu}{a}\right) \dddot{q} \, ,
\end{gather*}
where $a$, $b$, and $\mu$ are constants.
Replacing expression of $L$ and $Q$ into \eqref{eq:H14},
we conclude that Noether's nonconservative constants of motion
have the form
\begin{multline}
\label{eq:LC:ex3}
C(t,q,\dot{q},\ddot{q},\dddot{q})
= \frac{1}{2} \tau \left(\ddot{q}^2+a\dot{q}^2 + bq^2\right)
+(a\dot{q}-\dddot{q})(\xi_{0}-\dot{q}\tau)\\
+\ddot{q}(\xi_{1}-\ddot{q}\tau)-f + \Lambda\, ,
\end{multline}
where $f=\displaystyle \int(\xi-\dot{q}\tau)\left(\mu\dot{q}+\frac{\mu^2}{a^2}\ddot{q}
-\frac{2\mu}{a}\dddot{q}\right)dt$.
From the necessary and sufficient condition of invariance
\eqref{eq:H4}, one obtains that the exact symmetries ($\Lambda = 0$)
for the problem are given by
$\left(\tau,\xi_1,\xi_2,\sigma,\alpha_1,\alpha_2\right) = \left(c,0,0,0,0,0\right)$,
where $c$ is an arbitrary constant, and from \eqref{eq:LC:ex3} we conclude that
\begin{gather*}
\frac{1}{2} \left(\ddot{q}^2+a\dot{q}^2+bq^2\right)
-(a\dot{q}-\dddot{q})\dot{q}-\ddot{q}^2+\int\dot{q}\left(\mu\dot{q}
+\frac{\mu^2}{a^2}\ddot{q}-\frac{2\mu}{a}\dddot{q}\right)dt
\end{gather*}
is a conserved quantity for the nonconservative system.
This conclusion is nontrivial, and difficult to verify
directly from the definition of nonconservative constant
of motion (Definition~\ref{def:NCLs}).
\end{example}

The examples show how the previous nonconservative
results found in the literature are easily
covered by Theorem~\ref{th:TNNC}.
Theorem~\ref{th:TNNC} is, however, more general,
since it covers an arbitrary complex dynamical control system
of the form $\dot{q}(t)=\varphi\left(t,q(t),u(t)\right)$.
Moreover, Theorem~\ref{th:TNNC} introduces a new Hamiltonian
perspective to nonconservative Noether's theory.

%%%%%%%%%%%%%%%%%%%%%%%%%%%%%%%%%%%%%%%%%

\section*{Acknowledgments}

GF acknowledges the financial support of IPAD
(Portuguese Institute for Development);
DT the support from the Control Theory Group
(\textsf{cotg}) of the Centre for Research
in Optimization and Control of Aveiro University.

%%%%%%%%%%%%%%%%%%%%%%%%%%%%%%%%%%%%%%%%%

\end{document}